\documentclass[leqno,draft]{article}

\newtheorem{theorem}{Theorem}
\newtheorem{lemma}[theorem]{Lemma}
\newtheorem{proposition}[theorem]{Proposition}
\newtheorem{definition}[theorem]{Definition}
\newtheorem{corollary}[theorem]{Corollary}

\newcommand{\begintheorem}{\addtocounter{equation}{1}\begin{theorem}}
\newcommand{\beginlemma}{\addtocounter{equation}{1}\begin{lemma}}
\newcommand{\beginproposition}{\addtocounter{equation}{1}\begin{proposition}}
\newcommand{\begindefinition}{\addtocounter{equation}{1}\begin{definition}}
\newcommand{\begincorollary}{\addtocounter{equation}{1}\begin{corollary}}

\begin{document}

\title{An introduction to the geometry \\ of ultrametric spaces}

\author{Stephen Semmes \\
        Rice University}

\date{}

\maketitle

\begin{abstract}
Some examples and basic properties of ultrametric spaces are briefly
discussed.
\end{abstract}

\tableofcontents

\section{The ultrametric triangle inequality}
\label{ultrametric triangle inequality}
\setcounter{equation}{0}

         Let $(M, d(x, y))$ be a metric space.  Thus $M$ is a set, and
$d(x, y)$ is a nonnegative real-valued function defined for $x, y \in
M$ such that $d(x, y) = 0$ if and only if $x = y$,
\begin{equation}
         d(y, x) = d(x, y)
\end{equation}
for every $x, y \in M$, and
\begin{equation}
         d(x, z) \le d(x, y) + d(y, z)
\end{equation}
for every $x, y, z \in M$.  The latter condition is known as the
\emph{triangle inequality}, and if $d(x, y)$ satisfies the stronger condition
\begin{equation}
         d(x, z) \le \max(d(x, y), d(y, z))
\end{equation}
for every $x, y, z \in M$, then $d(x, y)$ is said to be an
\emph{ultrametric} on $M$.

         Equivalently, the ultrametric version of the triangle
inequality says that $d(x, y)$ and $d(y, z)$ cannot both be strictly
less than $d(x, z)$ for any $x, y, z \in M$.  In particular, the
standard metric on the real line does not have this property.

         The \emph{discrete metric} on any set $M$ is defined by
setting $d(x, y) = 1$ when $x \ne y$.  This is an ultrametric, and
there are also more interesting examples.

\section{$p$-Adic numbers}
\label{p-adic numbers}
\setcounter{equation}{0}

        Let $p$ be a prime number, such as $2$, $3$, $5$, or $7$, etc.
We would like to define the \emph{$p$-adic absolute value} $|\cdot
|_p$ on the set ${\bf Q}$ of rational numbers.  Let us begin by
defining $|\cdot |_p$ on the set ${\bf Z}$ of integers.

        If $x = 0$, then $|x|_p = 0$.  If $x \in {\bf Z}$ and $x \ne
0$, then
\begin{equation}
        x = w \, p^n
\end{equation}
for some $w, n \in {\bf Z}$ such that $w \ne 0$, $w$ is not divisible
by $p$, and $n \ge 0$.  In this case,
\begin{equation}
        |x|_p = p^{-n}.
\end{equation}
Note that
\begin{equation}
        |x|_p \le 1,
\end{equation}
since $n \ge 0$.

        Otherwise, if $x \in {\bf Q}$ and $x \ne 0$, then $x$ can be
expressed as
\begin{equation}
        x = \frac{u}{v} \, p^n
\end{equation}
for some $u, v, n \in {\bf Z}$ such that $u, v \ne 0$ and $u$, $v$ are
not divisible by $p$.  Again $|x|_p$ is defined to be $p^{-n}$, which
can be large or small depending on whether $n$ is negative or
positive.

        It is easy to see that
\begin{equation}
        |x + y|_p \le \max(|x|_p, |y|_p)
\end{equation}
for every $x, y \in {\bf Q}$.  If $x, y \in {\bf Z}$, then this
basically says that $x + y$ is divisible by $p^l$ when $x$ and $y$ are
both divisible by $p^l$.  Furthermore,
\begin{equation}
        |x \, y|_p = |x|_p \, |y|_p
\end{equation}
for every $x, y \in {\bf Q}$.

        The \emph{$p$-adic metric} is defined on ${\bf Q}$ by
\begin{equation}
        d_p(x, y) = |x - y|_p.
\end{equation}
This is an ultrametric on ${\bf Q}$, by the previous remarks.

        However, ${\bf Q}$ is not complete as a metric space with
respect to the $p$-adic metric.  In the same way that the real numbers
can be obtained by completing the rationals with respect to the
standard metric, the \emph{$p$-adic numbers} ${\bf Q}_p$ are obtained
by completing ${\bf Q}$ with respect to $d_p(x, y)$.

        As in the case of the real numbers, addition and
multiplication can be defined for $p$-adic numbers, and ${\bf Q}_p$ is
a field.  The $p$-adic absolute value and metric can be extended to
${\bf Q}_p$ as well, with properties like those just described on
${\bf Q}$.  The set ${\bf Z}_p$ of \emph{$p$-adic integers} consists
of $x \in {\bf Q}_p$ such that $|x|_p \le 1$, which one can show to be
the closure of ${\bf Z}$ in ${\bf Q}_p$.

\section{Spaces of sequences}
\label{sequence spaces}
\setcounter{equation}{0}

        Let $A$ be a set with at least two elements, and let
$\mathcal{A}$ be the set of sequences $a = \{a_i\}_{i = 1}^\infty$
such that $a_i \in A$ for each $i$.  If $0 < \rho < 1$ and $a, b \in
\mathcal{A}$, then put $d_\rho(a, b) = 0$ when $a = b$, and otherwise
\begin{equation}
        d_\rho(a, b) = \rho^n
\end{equation}
where $n$ is the largest nonnegative integer such that $a_i = b_i$
when $i \le n$.

        If $a, b, c \in \mathcal{A}$, then
\begin{equation}
        d_\rho(a, c) \le \max(d_\rho(a, b), d_\rho(b, c)).
\end{equation}
This basically follows from the fact that if $n$ is a nonnegative
integer such that
\begin{equation}
        a_i =b_i \hbox{ and } b_i = c_i \hbox{ when } i \le n,
\end{equation}
then $a_i = c_i$ for $i \le n$ too.  Thus $d_\rho(a, b)$ defines an
ultrametric on $\mathcal{A}$.

        Note that
\begin{equation}
        d_\rho(a, b) \le 1
\end{equation}
for each $a, b \in \mathcal{A}$ and $0 < \rho < 1$.  These
ultrametrics all determine the same topology on $\mathcal{A}$, which
is the product topology using the discrete topology on $A$.  If $A$
has only finitely many elements, then $\mathcal{A}$ is compact with
respect to this topology.  If we allow $\rho = 1$, then $d_\rho(a, b)$
reduces to the discrete metric on $\mathcal{A}$, and we get the
discrete topology on $\mathcal{A}$.  If $\rho > 1$, then $d_\rho(a,
b)$ is not even a metric on $\mathcal{A}$.

        For each $\alpha \in A$, let $T_\alpha : \mathcal{A} \to
\mathcal{A}$ be the mapping defined by $T_\alpha(a) = a'$, where
\begin{equation}
        a'_1 = \alpha \hbox{ and } a'_i = a_{i - 1} \hbox{ when } i \ge 2.
\end{equation}
It is easy to see that
\begin{equation}
        d_\rho(T_\alpha(a), T_\alpha(b)) = \rho \, d_\rho(a, b)
\end{equation}
for every $a, b \in \mathcal{A}$ and $0 < \rho < 1$.  Also,
\begin{equation}
        \mathcal{A} = \bigcup_{\alpha \in A} T_\alpha(\mathcal{A}).
\end{equation}
Thus $\mathcal{A}$ is the union of smaller copies of itself.

\section{Snowflake metrics}
\label{snowflakes}
\setcounter{equation}{0}

        Let $(M, d(x, y))$ be an ultrametric space.  For each $\tau >
0$, $d(x, y)^\tau$ is also an ultrametric on $M$, which determines the
same topology on $M$.  For example, if $d_\rho(a, b)$ is as in the
previous section, then
\begin{equation}
        d_\rho(a, b)^\tau = d_{\rho^\tau}(a, b).
\end{equation}
If $(M, d(x, y))$ is an ordinary metric space, then one can show that
$d(x, y)^\tau$ satisfies the triangle inequality and hence is a metric
on $M$ when $0 < \tau < 1$, but this does not always work when $\tau >
1$.  Suppose that $d(x, y)^\tau$ does satisfy the triangle inequality,
so that
\begin{equation}
        d(x, z) \le (d(x, y)^\tau + d(y, z)^\tau)^{1/\tau}
\end{equation}
for every $x, y, z \in M$.  This implies that
\begin{equation}
        d(x, z) \le 2^{1/\tau} \max(d(x, y), d(y,z))
\end{equation}
for every $x, y, z \in M$.  If $d(x, y)^\tau$ is a metric on $M$ for
every $\tau > 1$, then it follows that $d(x, y)$ is an ultrametric on
$M$, since $2^{1/\tau} \to 1$ as $\tau \to \infty$.

\section{Open sets}
\label{open sets}
\setcounter{equation}{0}

        Let $(M, d(x, y))$ be a metric space.  For each $x \in M$ and
$r > 0$, the open ball in $M$ with center $x$ and radius $r$ is defined by
\begin{equation}
        B(x, r) = \{y \in M : d(x, y) < r\}.
\end{equation}
Similarly, the closed ball in $M$ with center $x$ and radius $r$ is
defined by
\begin{equation}
        \overline{B}(x, r) = \{y \in M : d(x, y) \le r\}.
\end{equation}
A set $U \subseteq M$ is said to be \emph{open} if for every $x \in U$
there is an $r > 0$ such that
\begin{equation}
        B(x, r) \subseteq U.
\end{equation}
As usual, the empty set $\emptyset$ and $M$ itself are open subsets of
$M$, the intersection of finitely many open subsets of $M$ is an open
set, and the union of any family of open subsets of $M$ is an open set.

        For each $w \in M$ and $t > 0$, the open ball $B(w, t)$ is an
open set in $M$.  Indeed, let $x \in B(w, t)$ be given.  Thus $d(w, x)
< t$, and so $r = t - d(w, x) > 0$.  One can use the triangle
inequality to show that
\begin{equation}
        B(x, r) \subseteq B(w, t).
\end{equation}
If $d(\cdot, \cdot)$ is an ultrametric on $M$, then the same statement
holds with $r = t$.  Moreover,
\begin{equation}
        \overline{B}(x, t) \subseteq \overline{B}(w, t)
\end{equation}
for each $x \in \overline{B}(w, t)$.  Thus closed balls are also open
sets in ultrametric spaces.  Of course, this is not normally the case
in ordinary metric spaces, such as the real line with the standard metric.

        If $d(\cdot, \cdot)$ is an ultrametric on $M$, then the
complement of $B(w, t)$ in $M$ is also an open set, which is the same
as saying that $B(w, t)$ is a closed set in $M$.  Specifically, if $x
\in M \backslash B(w, t)$, then
\begin{equation}
        B(x, t) \subseteq M \backslash B(w, t).
\end{equation}
Equivalently, if $z \in M$ and $d(x, z) < r \le d(w, x)$, then
\begin{equation}
        d(w, z) \ge r.
\end{equation}
Otherwise, $d(w, z), d(x, z) < r$ imply that $d(w, x) < r$, a
contradiction.  Hence an ultrametric space $M$ with at least two
elements is not connected.  One can also use this to show that $M$ is
totally disconnected, in the sense that $M$ does not contain any
connected sets with at least two elements.  In particular, every
continuous path in an ultrametric space is constant.

\section{Completeness}
\label{completeness}
\setcounter{equation}{0}

        Let $(M, d(x, y))$ be a metric space.  Remember that a
sequence $\{x_j\}_{j = 1}^\infty$ of elements of $M$ is said to be a
\emph{Cauchy sequence} if for every $\epsilon > 0$ there is an $L \ge
1$ such that
\begin{equation}
        d(x_j, x_l) < \epsilon
\end{equation}
for every $j, l \ge L$.  If $\{x_j\}_{j = 1}^\infty$ is a Cauchy
sequence in $M$, then
\begin{equation}
        \lim_{j \to \infty} d(x_j, x_{j + 1}) = 0,
\end{equation}
which is to say that for every $\epsilon > 0$ there is an $L \ge 1$
such that
\begin{equation}
        d(x_j, x_{j + 1}) < \epsilon
\end{equation}
for each $j \ge L$.  The converse holds when $d(x, y)$ is an ultrametric, since
\begin{equation}
        d(x_j, x_l) \le \max(d(x_j, x_{j + 1}), \ldots, d(x_{l - 1}, x_l))
\end{equation}
when $j < l$.

        A metric space $M$ is said to be \emph{complete} if every
sequence of elements of $M$ converges to an element of $M$.  It is a
nice exercise to check that the spaces described in Section
\ref{sequence spaces} are complete.  Using the completeness of the
$p$-adic numbers ${\bf Q}_p$ and the preceding remarks, one can check
that an infinite series $\sum_{j = 1}^\infty x_j$ of $p$-adic numbers
converges in ${\bf Q}_p$ if and only if $\lim_{j \to \infty} x_j = 0$
in ${\bf Q}_p$.  If an infinite series $\sum_{j = 1}^\infty x_j$ of
real numbers converges, then $\lim_{j \to \infty} x_j = 0$ in ${\bf
R}$, but there are examples which show that the converse does not
hold.

\section{Binary sequences}
\label{binary sequences}
\setcounter{equation}{0}

        Let $\mathcal{B}$ be the set of all binary sequences, which is
to say the sequences $b = \{b_i\}_{i = 1}^\infty$ such that $b_i = 0$
or $1$ for each $i$.  Thus $\mathcal{B}$ is the same as the set
$\mathcal{A}$ associated to $A = \{0, 1\}$ as in Section \ref{sequence
spaces}.  Consider the mapping $\phi : \mathcal{B} \to {\bf R}$ that
sends each binary sequence to the real number with that binary
expansion.  Explicitly,
\begin{equation}
        \phi(b) = \sum_{i = 1}^\infty b_i \, 2^{-i}.
\end{equation}
It is well known that $\phi$ maps $\mathcal{B}$ onto the unit interval
$[0, 1]$, consisting of all real numbers $x$ such that $0 \le x \le
1$.  However, some real numbers have more than one binary expansion,
corresponding to binary sequences that are eventually constant.  This
is a continuous mapping with respect to the ultrametrics $d_\rho(a,
b)$ on $\mathcal{B}^*$ defined in Section \ref{sequence spaces} and
the standard metric on the real line.

        More precisely, let $d(a, b)$ be the ultrametric $d_\rho(a,
b)$ on $\mathcal{B}^*$ from Section \ref{sequence spaces} with $\rho =
1/2$.  Remember that the absolute value of a real number $x$ is
denoted $|x|$ and defined to be $x$ when $x \ge 0$ and $-x$ when $x
\le 0$, and that the standard metric on the real line ${\bf R}$ is
given by $|x - y|$.  For each integer $\ell$, one can check that
$\phi$ maps every closed ball in $\mathcal{B}^*$ with respect to $d(a,
b)$ of radius $2^\ell$ onto a closed interval in the real line with
length $2^\ell$.  In particular,
\begin{equation}
        |\phi(a) - \phi(b)| \le d(a, b)
\end{equation}
for every $a, b \in \mathcal{B}^*$.

\section{The Cantor set}
\label{cantor set}
\setcounter{equation}{0}

        If $r$, $t$ are real numbers with $r \le t$, then $[r, t]$ is
the closed interval in the real line consisting of $x \in {\bf R}$
such that $r \le x \le t$, and the length of this interval is $t - r$.
Let $E_0$ be the unit interval $[0, 1]$, and put
\begin{equation}
        E_1 = \Big[0, \frac{1}{3}\Big] \cup \Big[\frac{2}{3}, 1\Big].
\end{equation}
Continuing in this way, $E_n$ is the union of $2^n$ disjoint closed
intervals of length $3^{-n}$ for each positive integer $n$, and $E_{n
+ 1}$ is obtained from $E_n$ by removing the open middle third of each
of the $2^n$ intervals in $E_n$.  By construction, $E_{n + 1}
\subseteq E_n$ for each $n$, and $E = \bigcap_{n = 0}^\infty E_n$ is
known as the Cantor set.

        Consider the mapping $\psi : \mathcal{B} \to {\bf R}$ defined by
\begin{equation}
        \psi(b) = \sum_{i = 1}^\infty 2 \, b_i \, 3^{-i}.
\end{equation}
It is well known and not too difficult to show that $\psi$ is a
one-to-one mapping from $\mathcal{B}$ onto $E$.  Moreover, $\psi$ is a
homeomorphism with respect to the ultrametrics $d_\rho(a, b)$ on
$\mathcal{B}$ from Section \ref{sequence spaces}, $0 < \rho < 1$, and
the restriction of the standard metric on the real line to $E$.
For if $d(a, b) = d_\rho(a, b)$ on $\mathcal{B}$ with $\rho = 1/3$, then
\begin{equation}
        \frac{d(a, b)}{3} \le |\psi(a) - \psi(b)| \le d(a, b)
\end{equation}
for every $a, b \in \mathcal{B}$.  Indeed, if $a, b \in \mathcal{B}$
and $d(a, b) \le 3^{-n}$ for some integer $n \ge 0$, then $\psi(a)$
and $\psi(b)$ are in the same one of the $2^n$ intervals of length
$3^{-n}$ in $E_n$.  Hence the distance between $\psi(a)$ and $\psi(b)$
is at most $3^{-n}$.  If $d(a, b) = 3^{-n}$, then $\psi(a)$ and
$\psi(b)$ are in the two different subintervals of this interval that
are contained in $E_{n + 1}$, and the distance between $\psi(a)$ and
$\psi(b)$ is at least $3^{-n-1}$.

        The restriction of the standard metric on the real line to the
Cantor set is not an ultrametric.  However, the previous remarks show
that there is an ultrametric on the Cantor set which approximates the
standard metric, in the sense that each is bounded by a constant
multiple of the other.

\section{Doubly-infinite sequences}
\label{double sequences}
\setcounter{equation}{0}

        As in Section \ref{sequence spaces}, let $A$ be a set with at
least two elements, and choose an element $\alpha$ of $A$ to be a
basepoint.  Let $\mathcal{A}^*$ be the collection of doubly-infinite
sequences $a = \{a_i\}_{i = -\infty}^\infty$ of elements of $A$ for
which there is an integer $n$ such that $a_i = \alpha$ when $i \le n$.
If $0 < \rho < 1$, then $d_\rho(a, b)$ can be defined on
$\mathcal{A}^*$ as before by $d_\rho(a, b) = 0$ when $a = b$ and
$d_\rho(a, b) = \rho^n$ when $a \ne b$ and $n$ is the largest integer
such that $a_i = b_i$ for $i \le n$.  This is again an utrametric on
$\mathcal{A}^*$, which reduces to the discrete metric when $\rho = 1$.

        If we identify $\mathcal{A}$ with the set of $a \in
\mathcal{A}^*$ such that $a_i = \alpha$ when $i \le 0$, then this
definition of $d_\rho(a, b)$ agrees with the previous one on
$\mathcal{A}$.  The topology on $\mathcal{A}^*$ determined by
$d_\rho(a, b)$ is the same for each $\rho$, $0 < \rho < 1$, and
$d_\rho(a, b)$ is unbounded on $\mathcal{A}^*$.  If $A$ has only
finitely many elements, then $\mathcal{A}^*$ is locally compact with
respect to this topology.

        For each $a = \{a_i\}_{i = -\infty}^\infty \in \mathcal{A}^*$,
let $T(a)$ be the doubly-infinite sequence whose $i$th term is equal
to $a_{i - 1}$.  Thus $T(a) \in \mathcal{A}^*$ is the same as $a$, but
with the terms of the sequence shifted forward by one step.  This
defines a one-to-one mapping from $\mathcal{A}^*$ onto itself, which
satisfies
\begin{equation}
        d_\rho(T(a), T(b)) = \rho \, d_\rho(a, b)
\end{equation}
for every $a, b \in \mathcal{A}^*$ and $0 < \rho < 1$.  Note that
$T(a) = a$ exactly when $a_i = \alpha$ for each $i$.

\section{Binary sequences revisited}
\label{binary sequences, 2}
\setcounter{equation}{0}

        Let $\mathcal{B}^*$ be the set of doubly-infinite sequences $b
= \{b_i\}_{i = -\infty}^\infty$ such that $b_i = 0$ or $1$ for each
$i$ and there is an integer $n$ for which $b_i = 0$ when $i \le n$.
As in the previous section, we can identify the collection
$\mathcal{B}$ of binary sequences with the set of $b \in
\mathcal{B}^*$ that satisfy $b_i = 0$ when $i \le 0$.  We can also
extend the mapping $\phi$ from Section \ref{binary sequences} to
$\mathcal{B}^*$, by putting
\begin{equation}
        \phi(b) = \sum_{i = n}^\infty b_i \, 2^{-i}
\end{equation}
when $b_i = 0$ for each $i < n$.  Thus $\phi$ maps $\mathcal{B}^*$
onto the set of all nonnegative real numbers.  Two distinct elements
of $\mathcal{B}$ are sent to the same real number by $\phi$ if and
only if they agree up to some term, where one of the sequences is
equal to $0$ followed by all $1$'s, and the other is equal to $1$
followed by all $0$'s.  Let $d(a, b)$ be the ultrametric $d_\rho(a,
b)$ on $\mathcal{B}^*$ as in the previous section with $\rho = 1/2$.
As in Section \ref{binary sequences},
\begin{equation}
        |\phi(a) - \phi(b)| \le d(a, b)
\end{equation}
for every $a, b \in \mathcal{B}^*$.

\section{Nesting}
\label{nesting}
\setcounter{equation}{0}

        Let $(M, d(x, y))$ be an ultrametric space.  For every $x, y
\in M$ and $r, t > 0$, either
\begin{equation}
        B(x, r) \cap B(y, t) = \emptyset,
\end{equation}
or
\begin{equation}
        B(x, r) \subseteq B(y, t),
\end{equation}
or
\begin{equation}
        B(y, t) \subseteq B(x, r).
\end{equation}
More precisely, the first condition holds when
\begin{equation}
        d(x, y) > \max(r, t),
\end{equation}
the second condition holds when
\begin{equation}
        d(x, y) < t \hbox{ and } r \le t,
\end{equation}
and the third condition holds when
\begin{equation}
        d(x, y) < r \hbox{ and } t \le r.
\end{equation}
There are analogous statements for closed balls.  As usual, this does
not normally work in an ordinary metric space, like the real line with
the standard metric.  In an ultrametric space, it follows that
\begin{equation}
        B(x, r) \cup B(y, t)
\end{equation}
is a ball when
\begin{equation}
        B(x, r) \cap B(y, t) \ne \emptyset.
\end{equation}
This also holds in the real line, but not in the plane, for instance.

\section{Dyadic intervals}
\label{dyadic intervals}
\setcounter{equation}{0}

        If $r$, $t$ are real numbers with $r < t$, then the
half-open, half-closed interval $[r, t)$ consists of the real numbers
$x$ such that $r \le x < t$.  A \emph{dyadic interval} in the real
line is an interval of the form
\begin{equation}
        [i \, 2^l, (i + 1) \, 2^l),
\end{equation}
where $i$, $l$ are integers, although sometimes one considers closed
intervals of the same type instead.  Thus the length of a dyadic
interval is always an integer power of $2$.  If $I$, $I'$ are dyadic
intervals, then either $I \cap I' = \emptyset$, or $I \subseteq I'$,
or $I' \subseteq I$.  One should also include the possibility that
$I$, $I'$ are practically disjoint in the sense that $I \cap I'$
contains only a single point when one uses closed intervals.  The
structure of dyadic intervals basically corresponds to ultrametric
geometry, even if it may not be stated explicitly.  Dyadic intervals
and their relatives are often used in real analysis.

\section{More spaces of sequences}
\label{more sequence spaces}
\setcounter{equation}{0}

        Let $A_1, A_2, \ldots$ be a sequence of sets, each with at
least two elements.  As an extension of the situation described in
Section \ref{sequence spaces}, consider the space $\mathcal{A}$ of
sequences $a = \{a_i\}_{i = 1}^\infty$ such that $a_i \in A_i$ for
each $i$.  This is the same as the Cartesian product of the $A_i$'s.
Also let $\rho = \{\rho_i\}_{i = 0}^\infty$ be a strictly decreasing
sequence of positive real numbers with $\rho_0 = 1$.  For $a, b \in
\mathcal{A}$, put $d_\rho(a, b) = 0$ when $a = b$, and otherwise
\begin{equation}
        d_\rho(a, b) = \rho_n
\end{equation}
where $n$ is the largest nonnegative integer such that $a_i = b_i$ for
each $i \le n$.  This is equivalent to the earlier definition when
$\rho_n$ is the $n$th power of a fixed number.  As before, one can
check that $d_\rho(a, b)$ is an ultrametric on $\mathcal{A}$.  The
topology on $\mathcal{A}$ determined by this ultrametric is the
product topology using the discrete topology on each $A_i$ when
$\lim_{n \to \infty} \rho_n = 0$.  In this case, $\mathcal{A}$ is
compact if each $A_i$ has only finitely many elements.  There are
variants of this construction in which the geometry is less
homogeneous.

\end{document}